\documentclass[12pt,reqno]{amsart}

\usepackage{amsmath,amscd,amssymb,amsthm,mathtools}
\usepackage{mathdots}
\usepackage{bm}
\usepackage{yhmath}
\usepackage{todonotes} %add todo note
\usepackage[toc, page]{appendix}

\usepackage[linktocpage=true, backref=page]{hyperref}
\usepackage{mathrsfs}

\usepackage{derivative} %typeset differential operators

\usepackage{upgreek} % provide the upright Greek letters from the Euler or Adobe Symbol fonts as additional math symbols %e.g. \upmu, \uprho

\usepackage{microtype}

\usepackage{tikz,tikz-cd, color}
\usepackage[all]{xy}

\usepackage{wasysym, stackengine, makebox, graphicx}

\usepackage{adjustbox}
\usetikzlibrary{matrix}
\usetikzlibrary{decorations.pathmorphing}

\tikzset{
  symbol/.style={
    draw=none,
    every to/.append style={
      edge node={node [sloped, allow upside down, auto=false]{$#1$}}}
  }
}

\usepackage{stmaryrd}
\usepackage{enumerate}
\usepackage{paralist}   % for 'inparaenum' environment
\makeatletter
\newcommand{\myitem}[1]{%
\item[(#1)]\protected@edef\@currentlabel{#1}%
}
\makeatother
\usepackage{xcolor}
\usepackage{booktabs}
\usepackage{float}
\mathchardef\mhyphen="2D % Define a "math hyphen"

\usepackage{fullpage}

\newtheorem{thm}{Theorem}
\newtheorem*{thm*}{Theorem}

\newtheorem{corollary}[thm]{Corollary}
\newtheorem{thm&defn}[thm]{Theorem \& Definition}

\newtheorem{theorem}[thm]{Theorem}

\newtheorem{proposition}[thm]{Proposition}

\newtheorem{conjecture}[thm]{Conjecture}

\theoremstyle{definition}

\newtheorem{definition}[thm]{Definition}
\newtheorem*{notation*}{Notation}
\newtheorem{notation}[thm]{Notation}

\newtheorem*{convention*}{Convention}
\newtheorem{example}[thm]{Example}

\newtheorem*{claim*}{Claim}

\theoremstyle{remark}

\newtheorem*{rem*}{Remark}
\newtheorem{remark}[thm]{Remark}

\numberwithin{equation}{section}
\numberwithin{thm}{section}
\newcommand{\bZ}{\mathbb Z}
\newcommand{\bN}{\mathbb N}
\newcommand{\bQ}{\mathbb Q}

\newcommand{\bC}{\mathbb C}

\newcommand{\bP}{\mathbb{P}}

\newcommand{\bF}{\mathbb{F}}

\newcommand{\cB}{\mathcal{B}}

\newcommand{\cO}{\mathscr{O}}
\newcommand{\cS}{\mathcal{S}} % deformation family
\newcommand{\cW}{\mathcal{W}} %bundle %family
\newcommand{\cX}{\mathcal{X}}
\newcommand{\cY}{\mathcal{Y}}

\newcommand{\fX}{\mathfrak{X}}

\newcommand{\fz}{\mathfrak{z}}

\newcommand{\ff}{\mathfrak{f}}

\newcommand{\ZGWp}{\mathrm{Z}_\GW^\pp}
\newcommand{\ZPT}{\mathrm{Z}_\PT}

\newcommand{\tmid}{\,\middle\vert\,}

\newcommand{\pp}{\prime} % disconnected

\newcommand{\sC}{\mathsf{c}}

\newcommand{\tsK}{\widetilde{\mathsf{K}}}

\newcommand{\sZ}{\mathsf{Z}}

\newcommand{\oA}{\overline{A}}
\newcommand{\oM}{\overline{\mathcal{M}}} %moduli
\newcommand{\oMp}{\overline{\mathcal{M}}^\pp} %moduli (possibly disconnected)
\newcommand{\oY}{\bar{Y}}
\newcommand{\halpha}{\widehat{\alpha}}
\newcommand{\hell}{\widehat{\ell}}

\newcommand{\PSR}[1]{\llbracket #1\rrbracket}
\newcommand{\Laurent}[1]{(\hspace{-0.25em}(#1)\hspace{-0.25em})}

\newcommand{\Xres}{Y} %resoltion space in a transition
\newcommand{\Xsing}{\bar{Y}} %singular space in a transition
\newcommand{\Xsm}{X} %smoothing space in a transition
\newcommand{\crpcon}{\phi} %crpeant contraction in a transition

\makeatletter % pmatrix[some constant]
\renewcommand*\env@matrix[1][\arraystretch]{%
  \edef\arraystretch{#1}%
  \hskip -\arraycolsep
  \let\@ifnextchar\new@ifnextchar
  \array{*\c@MaxMatrixCols c}}
\makeatother

\newcommand{\Def}{\mathrm{Def}}

\newcommand{\ev}{\mathrm{ev}}
\newcommand{\loc}{\mathrm{loc}}
\newcommand{\vir}{\mathrm{vir}}
\newcommand{\GW}{\mathrm{GW}}
\newcommand{\PT}{\mathrm{PT}}
\newcommand{\vacuum}{|0\rangle}
\newcommand{\tp}{\tau} %pt descendent insertion
\newcommand{\pt}{\mathrm{pt}}

\newcommand{\chitop}{\chi_\mathrm{top}} %topological Euler characteristic

\newcommand{\typeIIno}{\Xres \rightarrow \Xsing \rightsquigarrow \Xsm}
\newcommand{\typeII}{\Xres \xlongrightarrow{\crpcon} \Xsing \rightsquigarrow \Xsm} %Type II extremal transition

\newcommand{\Xloc}{X_\loc}
\newcommand{\Yloc}{Y_\loc}

\newcommand{\I}{\mathrm{I}}
\newcommand{\II}{\mathrm{I\!I}}

\DeclareMathOperator{\Aut}{Aut}
\DeclareMathOperator{\NE}{NE} %Mori cone
\DeclareMathOperator{\ch}{ch} %chern character

\DeclareMathOperator{\Bl}{Bl} %blow-up

\DeclareMathOperator{\Hilb}{Hilb}

\newcommand{\ignore}[1]{}

\title{Remarks on GW/PT under del Pezzo transitions}

\author{Shuang-Yen Lee}
\address[]{Department of Mathematics, National Taiwan University, Taipei 10617, Taiwan}
\email{d10221004@ntu.edu.tw}

\author{Chin-Lung Wang}
\address[]{Department of Mathematics and Taida Institute for Mathematical Sciences (TIMS), National Taiwan University, Taipei 10617, Taiwan}
\email{dragon@math.ntu.edu.tw}

\author{Sz-Sheng Wang}
\address[]{Department of Applied Mathematics, National Yang Ming Chiao Tung University, Hsinchu 30010, Taiwan}
\email{sswangtw@math.nctu.edu.tw}


\begin{document}

\begin{abstract} 
A projective threefold transition $Y \xlongrightarrow{\phi} \bar{Y} \rightsquigarrow X$ is del Pezzo if $\phi$ contracts a smooth del Pezzo surface to a point. We show that the GW/PT correspondence holds on $Y$ implies that it holds on $X$. In particular, a hypersurface of degree $6$ in $\mathbb{P} (3, 2, 1, 1, 1)$ gives a new example to the correspondence. The main tools are (i) the double point degeneration constructed in \cite{LWW25} and (ii) deformations of del Pezzo surfaces into toric surfaces (Proposition~\ref{prop;def_to_toric}). Applications of the degeneration formulas in GW and PT then reduce the problem to known cases.  
\end{abstract}

\maketitle

%%%%%%%%%%%%%%%%%%%%%%%%%%%%%%%%%%%%%%%%
%%%%%%%%%%%%%%%%%%%%%%%%%%%%%%%%%%%%%%%%
\section{Introduction}
%%%%%%%%%%%%%%%%%%%%%%%%%%%%%%%%%%%%%%%%
%%%%%%%%%%%%%%%%%%%%%%%%%%%%%%%%%%%%%%%%

Two smooth projective threefolds $\Xsm$ and $\Xres$ are related by a \emph{geometric transition} if there exists a crepant contraction $\crpcon \colon \Xres \to \Xsing$ followed by a smoothing $\fX \to \Delta$ of $\Xsing = \fX_0$ with the fiber $X = \fX_t$ for some $t \neq 0$. We write $\Xres \searrow \Xsm$  or $\typeII$ for this process.

When the singular variety $\Xsing$ has only ordinary double points, this is known as a conifold transition. Comparison of geometric invariants under conifold transitions has been intensively studied in the literature. For example, the Gromov--Witten/Pandharipande--Thomas (GW/PT) correspondence (Conjecture~\ref{GWPTconj}) holds on $Y$ implies that it holds on $X$ \cite{LW25}.

Here we are interested in more general transitions. A transition $\typeII$ is a \emph{del Pezzo transition of degree $d$} if the $\crpcon$ exceptional set $E$ is a smooth divisor with $E^3 = d$ and $\crpcon (E)$ is a point. In particular $E \cong S_d$ is a smooth del Pezzo surface of degree $d$. The purpose of this note is to compare certain geometric invariants under del Pezzo transitions.  
 
The main tool used is the double point degenerations we have found in \cite[\S 4]{LWW25} (see \eqref{eqn;ss_Y} and \eqref{eqn;ss_X}). In simple terms, to conclude that the semistable model of the smoothing $\fX \to \Delta$ is of double point type, we control the base change degree for $\fX \to \Delta$ so that we may first perform the base change and then perform only one weighted blow-up to achieve the semistable model (see \eqref{eqn;ss_X_diag}). The double point degenerations then allow us to employ various existing degeneration formulas to study the comparison by reducing the problem to its local models $Y_d \searrow X_d$ for $d \in \{1, 2, 3, 4, 5, 6\I, 6\II, 7, 8\}$. Here $X_d$ (resp.\ $S_d$) is a smooth del Pezzo threefold (resp.\ surface) of degree $d$ , $Y_d = \bP_{S_d}(K_{S_d} \oplus \cO)$. Note that for $d = 6$ there are indeed two possible smoothings (see Example \ref{ex:dPtrans}).
 
As a first application, we prove in Theorem \ref{thm;top} that if $\varphi$ is a group homomorphism from the complex cobordism group to a group, then 
\begin{align*}
\varphi(\Xres) - \varphi (\Xsm) = \varphi(\Xres_d) - \varphi (\Xsm_d). 
\end{align*}
In particular, up to complex cobordism, a del Pezzo transition is equivalent to its local model. 

Secondly, we study the descendent GW/PT correspondence (Conjecture~\ref{GWPTconj}) under del Pezzo transitions. We show in Theorem~\ref{thm;GWPT} that if $\Xres$ satisfies Conjecture~\ref{GWPTconj}, then so does $\Xsm$ for descendent insertions \eqref{eqn;ins-total}. To apply degeneration formulas for GW and PT-invariants we need to prove that the relative descendent GW/PT correspondence (Conjecture~\ref{relGWPTconj}) holds for the pairs $(Y_d, E)$ and $(Y_d, H)$, where $E$ (resp.\ $H$) is the zero (resp.\ infinity) section of $Y_d$ (see Theorem~\ref{thm;rel_loc_GWPTdPs}). The trick is to show that that every del Pezzo surface is deformation equivalent to a toric surface (Proposition~\ref{prop;def_to_toric}) and therefore $Y_d$ is deformation equivalent to a toric threefold whose GW/PT correspondence is known by \cite{PP14}. 

As a byproduct, we also prove Conjecture~\ref{GWPTconj} for certain (weak) Fano threefolds (Theorem~\ref{thm;GWPTcorrp_def_toric}). In particular we obtain new cases which are not contained in \cite{LW25,PP14,PP17}, e.g.~the del Pezzo threefold of degree one $X_1$, which is a hypersurface of degree $6$ in $\bP (3, 2, 1, 1, 1)$ (see Remarks~\ref{rmk;exa} and \ref{rmk;exaX1}).

%%%%%%%%%%%%%%%%%%%%
\subsection*{Acknowledgments}

C.-L.~is supported by NSTC, Taiwan with grant number NSTC 113-2115-M-002-004-MY3, and by the Core Research Group of National Taiwan University. S.-Y.~is supported by the PhD fund in C.-L.'s NSTC grant above. S.-S.~is supported by NSTC, Taiwan with grant numbers NSTC 111-2115-M-A49-019-MY3 and 114-2115-M-A49-007-MY3.
%the National Science and Technology Council
%%%%%%%%%%%%%%%%%%%%

%%%%%%%%%%%%%%%%%%%%%%%%%%%%%%%%%%%%%%%%
%%%%%%%%%%%%%%%%%%%%%%%%%%%%%%%%%%%%%%%%
\section{Double Point Degenerations}\label{sec;dP}
%%%%%%%%%%%%%%%%%%%%%%%%%%%%%%%%%%%%%%%%
%%%%%%%%%%%%%%%%%%%%%%%%%%%%%%%%%%%%%%%%

%%%%%%%%%%%%%%%%%%%%%%%%%%%%%%%%%%%%%%%%
\subsection{Del Pezzo Transitions}\label{subsusec;dP_transition}
%%%%%%%%%%%%%%%%%%%%%%%%%%%%%%%%%%%%%%%%

We briefly review the basics of del Pezzo transitions (see \cite[\S 1 \& \S 4]{LWW25} and references therein for more details) and set notations for the rest of this paper.

Given a del Pezzo transition $\typeII$ of degree $d$, we know that the exceptional divisor $E$ is a smooth del Pezzo surface of degree $d$ (cf.~\cite[Proposition 2.13]{Reid80}). Moreover, $E$ is not isomorphic to $\bP^2$ or the Hirzebruch surface $\bF_1$ because the deformation space $\Def(\Xsing, p)$ must contain a smoothing component (cf.~\cite[Remark 1.8]{LWW25}), where $\crpcon (E) = \{p\}$. In particular we have $1 \leq d \leq 8$.

We shall recall the construction of a (standard) local model $Y_d \searrow X_d$ of del Pezzo transitions degree of $d$ in the following example.

\begin{example}\label{ex:dPtrans}
For $1 \leq d \leq 8$, let $S_d$ denote a smooth del Pezzo surface of degree $d$ which is isomorphic to the blow-up of $\bP^1 \times \bP^1$ at $8 - d$ points. 

Let us denote by $\alpha$ the weight
\begin{align}\label{eqn;wt}
    \alpha =
    \begin{cases}
       (3, 2, 1, 1) & \text{if } d = 1,\\
       (2, 1, 1, 1) & \text{if } d = 2,\\
       (1, \ldots, 1) & \text{if } 3 \leq d \leq 8,
    \end{cases}
\end{align}
where the last sequence of $1$ is repeated $d + 1$ times. Then we have the anti-canonical embedding $S_d \hookrightarrow \bP (\alpha)$. Let $\oY_d$ be the projective cone over $S_d$ with vertex $p$ in the weighted projective space $\bP (\alpha, 1)$ and 
\begin{align}\label{eqn;def_Yd}
    Y_d = \bP_{S_d}(K_{S_d} \oplus \cO).
\end{align}
It is immediate that $Y_d$ is the weighted blow-up of $\oY_d$ at the point $p = [0: \dots :0:1] \in \bP (\alpha, 1)$ with the weight $\alpha$. Therefore $S_d$ is the exceptional divisor of $Y_d \to \oY_d$. Note that $Y_d \to \oY_d$ is the restriction of the weighted blow-up $\bP_{\bP (\alpha)} (\cO(-1) \oplus \cO) \to \bP (\alpha, 1)$ at $p$ with the weight $(\alpha, 1)$. To summarize, we have
\begin{equation*}
    \begin{tikzcd}
       Y_d \ar[r,hook] \ar[d] & \bP_{\bP (\alpha)} (\cO(-1) \oplus \cO)  \ar[d]  \\
        \oY_d \ar[r,hook] & \bP (\alpha, 1).
    \end{tikzcd}
\end{equation*}

We denote by $(X_d, \cO_{X_d} (1))$ a smooth del Pezzo threefold of degree $d$. By the classification of Fujita and Iskovskikh (see \cite{bookAGV}, \cite[Appendix A]{LWW25} and the references therein), we have the anti-canonical embedding $X_d \hookrightarrow \bP(\alpha, 1)$. Moreover, $X_d$ is one of the following:
\begin{enumerate}
    \item $d = 1$ and $X_1$ is a hypersurface of degree $6$ in $\bP (3, 2, 1, 1, 1)$.

    \item $d = 2$ and $X_2$ is a double cover $X \to \bP^3$ ramified along a surface in $\bP^3$ of degree $4$.
    %quartic hypersurface in $\bP (2, 1, 1, 1, 1)$.

    \item $d = 3$ and $X_3 \hookrightarrow \bP^{4}$ is a hypersurface of degree $3$.

    \item $d = 4$ and $X_4 \hookrightarrow \bP^{5}$ is a complete intersection of two quadrics.

    \item $d = 5$ and $X_5 \hookrightarrow \bP^{6}$ is a linear section of Pl\"ucker embedding of $\mathrm{Gr}(2,5) \subseteq \bP^9$ by a codimension $3$ subspace.
    
    \myitem{6I}\label{thm;dP3claf_6I} $d = 6$ and $X_{6\I} \hookrightarrow \bP^{7}$ is a hypersurface of bidegree $(1,1)$ in $\bP^2 \times \bP^2$ , and $\bP^2 \times \bP^2 \hookrightarrow \bP^8$ by Segre embedding.
    
    \myitem{6I\!I}\label{thm;dP3claf_6II} $d = 6$ and $X_{6\II} = \bP^1 \times \bP^1 \times \bP^1 \hookrightarrow \bP^7$ by Segre embedding. 

    \setcounter{enumi}{6}
    
    \item $d = 7$ and $X_7 \hookrightarrow \bP^{8}$ is a blow-up of $\bP^3$ at a point.

    \item $d = 8$ and $X_8 = \bP^{3} \hookrightarrow \bP^{9}$ with $\cO_X (1) = \cO_{\bP^3} (2)$.
\end{enumerate}
Then $X_d$ is a smoothing of the singular del Pezzo threefold $\oY_d$ and therefore we get the desired del Pezzo transition $Y_d \searrow X_d$ for $d \in \{1, 2, 3, 4, 5, 6\I, 6\II, 7, 8\}$. Here we adopt the convention that $Y_{6\I} = Y_{6\II} \coloneqq Y_6$.
\end{example}

Let $\fX \to \Delta$ be the corresponding smoothing of the del Pezzro transition $\typeII$.
Note that $\crpcon$ is the weighted blow-up at the unique singularity $p$ of $\Xsing$ with weight $\alpha$ \eqref{eqn;wt} (cf.~\cite[Theorem 2.11]{Reid80}).

To use the local model $Y_d \searrow X_d$ to study $\Xres \searrow \Xsm$, we need two double point degenerations. Let us review the construction of such degenerations in \cite[\S 4.1]{LWW25}.

The \emph{K\"ahler degeneration} $\cY \coloneqq \Bl_{E \times \{0\}}( \Xres \times \Delta) \to \Delta$ is the deformation to the normal cone. Since $E$ has codimension one in $\Xres$ and $E|_E = K_E$, the special fiber $\cY_0 = \Xres \cup Y_d$ is a simple normal crossing divisor with $Y_d = \bP_E (K_E \oplus \cO)$. The intersection $E = \Xres \cap Y_d$ is understood as the infinity divisor (or relative hyperplane section) of $Y_d \to E$. We also denote this degeneration by 
\begin{align}\label{eqn;ss_Y}
    Y \rightsquigarrow Y \cup_{E} Y_d.
\end{align}

The \emph{complex degeneration} is $\cX \to \Delta$ is the semistable reduction of the smoothing $\fX \to \Delta$. Set $n_d = 6, 4, 3$ for $d = 1,2,3$ respectively, and $n_d = 2$ for $d \geq 4$. It is obtained by a degree $n_d$ base change $\cX' \to \Delta$ allowed by the weighted blow-up at $p \in \cX'$ with weight $(\alpha, 1)$: 
\begin{equation}\label{eqn;ss_X_diag}
    \begin{tikzcd}
        \cX \ar[r] \ar[d]& \fX' \ar[r] \ar[d] \ar[dr, phantom, "\square"]& \fX \ar[d]\\
        \Delta \ar[r, equal]& \Delta \ar[r] & \Delta.
    \end{tikzcd}
\end{equation}
For the control of the base change degree $n_d$, see the proof of Proposition 1.10 in \cite{LWW25}. The special fiber $\cX_0 = \Xres \cup X_d$ is a simple normal crossing divisor with $X_d$ being a smooth del Pezzo threefold of degree $d$. The intersection $E = \Xres \cap X_d$ in $X_d$ is now understood as a general member of the linear system $|{-K_{X_d}}|$, and the normal bundle of the intersection $E$ in $Y$ is $K_E$. In particular, we have
\begin{align}\label{eqn;nomal_bd}
    N_{E/ Y} \otimes N_{E/X} \cong \cO(K_E) \otimes \cO(- K_E) \cong \cO_E.
\end{align}
We also denote this degeneration by 
\begin{align}\label{eqn;ss_X}
    X \rightsquigarrow Y \cup_{E} X_d.
\end{align}

Notice that the local model $Y_d \searrow X_d$ appears in the special fibers $\cX_0$ and $\cY_0$. For $d = 6$, $X_{6\I}$ and $X_{6\II}$ are distinguished by the irreducible component of $\Def(\Xsing, p)$ that contains the image of the holomorphic map $\Delta \to \Def(\Xsing, p)$ induced by the smoothing $\fX \to \Delta$.

%%%%%%%%%%%%%%%%%%%%%%%%%%%%%%%%%%%%%%%%
\subsection{Topology: Chern Numbers}
%%%%%%%%%%%%%%%%%%%%%%%%%%%%%%%%%%%%%%%%

Applying the double point degeneration \eqref{eqn;ss_X} to del Pezzo transitions of degree $d$, we can identify them with the local model $Y_d \searrow X_d$ in the complex cobordism ring $\Omega^U_\ast$. Recall that $\Omega^U_\ast$ is generated by all stable almost complex manifolds, and it is a polynomial ring over $\bZ$. Two manifolds determine the same element in $\Omega^U_\ast$ if and only if their Chern numbers coincide.

\begin{theorem}\label{thm;top}
Given a del Pezzo extremal transition $\typeII$ of degree $d$ , we let $Y_d$ and $X_d$ be as in Example \ref{ex:dPtrans}. If $\varphi \colon \Omega^U_6 \to \Lambda$ is a group homomorphism, then 
\begin{align*}%\label{eqn;top}
    \varphi(\Xres) - \varphi (\Xsm) = \varphi(\Xres_d) - \varphi (\Xsm_d) 
\end{align*}
for $d \in \{1, 2, 3, 4, 5, 6\I, 6\II, 7, 8\}$.
\end{theorem}

\begin{proof}
We will apply algebraic cobordism theory to prove this theorem, see \cite{LM07, LP09} for further details and references. Indeed, we use the double point cobordism ring $\omega_\ast(\bC)$ of Levine-Pandharipande \cite[Theorem 1]{LP09} to give a simple description of the difference $[X] - [Y]$.

By \eqref{eqn;nomal_bd}, \eqref{eqn;ss_X} and double point relations \cite[Definition 0.1]{LP09}, we get 
\begin{align}\label{eqn;doublpt_rel}
    [X] - [Y] - [X_d] + [\bP (N_{E / Y} \oplus \cO)] = 0
\end{align}
in $\omega_3 (\bC)$. According to \eqref{eqn;def_Yd} and that $\crpcon$ is crepant, it follows that $\bP (N_{E / Y} \oplus \cO) = Y_d$ by the adjunction formula. Then the theorem follows from \eqref{eqn;doublpt_rel}, the isomorphism 
\[
    \omega_\ast (\bC) \cong \Omega^U_{2 \ast}
\]
(see \cite[Lemma 4.3.1, Theorem 4.3.7]{LM07}) and that $\varphi$ is a group homomorphism.
\end{proof}

We give some examples of group homomorphims $\varphi$.

\begin{example}
If the group $\Lambda$ is $\Omega_6^U$, then the identity map $\varphi = \mathrm{id}$ is a trivial example, and we get \eqref{eqn;doublpt_rel}.

Let $\sZ(M, q)$ be the partition function for degree $0$ Donaldson--Thomas invariants on a smooth projective threefold $M$. It gives another example of group homomorphisms. In fact, by the degeneration formula in Donaldson--Thomas theory, a group homomorphism (see \cite[\S 13]{LP09}) 
\begin{align*}
    \varphi \colon \Omega_6^U \cong \omega_3 (\bC) \to \bQ \PSR{q}^\ast
\end{align*}
is defined by $\varphi (M) \coloneqq \sZ(M, q)$, where $\bQ \PSR{q}^\ast \subseteq\bQ \PSR{q}$ is the multiplicative group of power series with constant term $1$.
\end{example}

\begin{example}
Let $R$ be a commutative ring containing $\bQ$. We recall the Hirzebruch $R$-genus $\varphi$ (see \cite[\S 1]{Hirzebruch95}). By definition, it is a ring homomorphism $\varphi \colon \Omega_\ast^U \otimes \bQ \to R$, which depends only on Chern numbers. To each series of the form $Q(x) = 1 + a_1 x + a_2 x^2 + a_3 x^3 + \cdots \in R\PSR{x}$ there corresponds the Hirzebruch $R$-genus $\varphi_Q (M) \coloneqq \int_M \prod_{i} Q (x_i)$, where $x_i$'s are the Chern roots of a stable almost complex manifold $M$. If $M$ has complex dimension three, then
\begin{align}\label{eqn;Rgenus}
    \varphi_Q (M) = (a_1^3 - a_1a_2 + a_3) c_3 + (a_1a_2 - 3a_3) c_1c_2 + a_3 c_1^3,
\end{align}
where $c_i = c_i (M)$ is the Chern number.

Since $\varphi_Q$ is also a group homomorphism, we can apply Theorem \ref{thm;top} to compute the difference of the $R$-genus of a del Pezzo transition $\Xres \searrow \Xsm$ of degree $d$. It suffices to compute the Chern numbers of the local model $Y_d \searrow X_d$. By \eqref{eqn;def_Yd} and Euler sequence for the projective bundle $\pi \colon Y_d \to S_d$,
\[
    0 \to \cO_{\Xres_d} \to \pi^\ast(K_{S_d} \oplus \cO) \otimes \cO_{\Xres_D}(1) \to T_{\Xres_d} \to \pi^\ast T_{S_d} \to 0,
\]
it is easily seen that 
$c_1(Y_d)^3 =  8d$, $c_1(Y_d) c_2(Y_d) = 24$ and
\begin{align*}
    c_3 (Y_d)= 2(12-d) = 2 c_2(S_d),
\end{align*}
where $S_d$ is a smooth del Pezzo surface of degree $d$. On the other hand, since $X_d$ is a del Pezzo threefold of degree $d$, we have that $c_1(X_d)^3 = 8 d$. From standard arguments using the Riemann–Roch, Serre duality, and Kodaira vanishing  \cite[Corollary 2.1.14]{bookAGV}, we find that
\[
    d + 2 = \chi (\cO (1)) = d + \frac{1}{12} c_1(X_d) c_2(X_d),
\]
i.e., $c_1(X_d) c_2(X_d) = 24$. The top Chern number of $X_d$ is given in Table \ref{tab;topXd} by the classification of del Pezzo threefolds, see for example \cite[Appendix A]{LWW25} and references therein. Set $\Delta\chitop \coloneqq \chitop (Y_d) - \chitop(X_d)$, where $\chitop(-)$ is the topological Euler number. We also list the topological difference $\Delta\chitop$ of $\Xres \searrow \Xsm$ in Table \ref{tab;topXd} (cf.~\cite[Remark 1.12]{LWW25}). Therefore 
\begin{align*}
    \varphi_Q(\Xres) - \varphi_Q (\Xsm) &= (a_1^3 - a_1a_2 + a_3) (2 \chitop(S_d) - \chitop(X_d)) \\
    &= (a_1^3 - a_1a_2 + a_3) \Delta\chitop
\end{align*}
by \eqref{eqn;Rgenus} and $\chitop(Y_d) = \chitop(\bP^1)\chitop(S_d)$. As a byproduct of \ref{eqn;doublpt_rel} and that $X_d$ and $Y_d$ have same $c_1c_2 = 24$, we also obtain that $c_1(\Xres) c_2(\Xres) = c_1(\Xsm) c_2(\Xsm)$.

\begin{table}[H]
    \centering
    \begin{tabular}{c||ccccccccc}
        \toprule
        $d$ & $1$ & $2$ & $3$ & $4$ & $5$ & $6 \I$ & $6 \II$ & $7$ & $8$ \\
        \midrule
        $c_3 (X_d)$ & $- 38$ & $-16$ & $-6$ & $0$ & $4$ & $6$ & $8$ & $6$ & $4$ \\
        $\Delta \chitop$ & $60$ & $36$ & $24$ & $16$ & $10$ & $6$ & $4$ & $4$ & $4$ \\ 
        \bottomrule
    \end{tabular}
    \caption{Topological numbers of $Y_d \searrow X_d$.}
  \label{tab;topXd}
\end{table}
\end{example}

%%%%%%%%%%%%%%%%%%%%%%%%%%%%%%%%%%%%%%%%
%%%%%%%%%%%%%%%%%%%%%%%%%%%%%%%%%%%%%%%%
\section{Quantum: GW/PT Correspondence}\label{sec;GWPT}
%%%%%%%%%%%%%%%%%%%%%%%%%%%%%%%%%%%%%%%%
%%%%%%%%%%%%%%%%%%%%%%%%%%%%%%%%%%%%%%%%
In this section, we will use the double point degenerations \eqref{eqn;ss_Y} and \eqref{eqn;ss_X} to relate descendent GW/PT correspondence under del Pezzo transitions (Theorem \ref{thm;GWPT}). We also prove correspondences for certain (weak) Fano threefolds (Theorems \ref{thm;GWPTcorrp_def_toric} and \ref{thm;rel_loc_GWPTdPs}). As \S \ref{subsusec;abs_theories} and \ref{subsusec;rel_theories} are primarily to review the basics of GW and PT theories and to set notations, the exposition is condensed. See \cite{PP14, PP17, LW25} and references therein for more information.

%%%%%%%%%%%%%%%%%%%%%%%%%%%%%%%%%%%%%%%%
\subsection{Absolute Theories}\label{subsusec;abs_theories}
%%%%%%%%%%%%%%%%%%%%%%%%%%%%%%%%%%%%%%%%

Let $M$ be a smooth projective threefold. Fix a curve class $0 \neq \beta\in \NE(M)$, integers $r\in \mathbb{Z}_{\ge 0}$ and $g\in \mathbb{Z}$. Set $\sC_\beta \coloneqq ( c_1 (T_V), \beta)$.

Let $\oM^\pp_{g, r} (M, \beta)$ denote be the moduli space of $r$-marked genus $g$ degree $\beta$ stable maps $C \to M$, where the stable map is required to have positive degree on each connected component of the (possibly disconnected) domain $C$. The moduli space $\oM^\pp_{g, r} (M, \beta)$ is equipped with a virtual fundamental class and its virtual dimension is $\sC_\beta + r$. 

\begin{definition}
Let $\psi_j$ be the first Chern class of cotangent line bundle associated to the $j$-th marked point for $j = 1$, $\dots$, $r$. Then the disconnected descendent GW-invariant is defined as
\begin{align*}
    \langle \tau_{\alpha_1-1} (\gamma_1) \cdots \tau_{\alpha_r-1} (\gamma_r) \rangle^\pp_{g, \beta} = \int_{[\oM^\pp_{g, r} (M, \beta)]^\vir} \prod_{j = 1}^r 
    \psi_j^{\alpha_j-1} \cup \ev_j^\ast (\gamma_j)
\end{align*}
where $\ev_j$ is the evaluation map given by the $j$-th marked point and $\gamma_j \in H^\ast (M, \bQ)$.
We define the following associated partition function
\begin{align} \label{eqn;GW_asso_fun}
    \ZGWp \left(M; u \tmid  \prod_{j = 1}^r \tau_{\alpha_j-1} (\gamma_j) \right)_\beta = \sum_{g \in \bZ} \left\langle \prod_{j = 1}^r \tau_{\alpha_j-1} (\gamma_j) \right\rangle^\pp_{g, \beta} u^{2g - 2}.
\end{align}
\end{definition}

Note that $\oM^\pp_{g, r} (M, \beta)$ is empty for $g$ sufficiently negative. Therefore \eqref{eqn;GW_asso_fun} is a Laurent series in $\bQ\Laurent{u}$.

To define PT-invariants, we consider the moduli space of stable pairs. A \emph{stable pair} 
\begin{align*}
    (F, s\colon \cO_M \to F)
\end{align*} 
on $M$ consists of a pure one-dimensional sheaf $F$ on $M$ and a section $s$ with zero-dimensional cokernel. 
Given $n \in \bZ$, let $P_n (M, \beta)$ be the moduli space of stable pairs with $\ch_2 (F) = \beta$ and $\chi (F) = n$. Then $P_n (M, \beta)$ is fine and projective, and it admits a virtual fundamental class of virtual dimension $\sC_\beta$. %\cite[Theorem 2.14]{PT09}
Let $\bF$ be the universal sheaf of $P_n (M, \beta)$. Consider the $k$-th descendent insertion 
\begin{align*}
    \tp_k (\gamma) \coloneqq \pi_{P \ast} (\pi_M^\ast (\gamma) \cdot \ch_{2 + k} (\bF)) \in H^\ast (P_n (M, \beta), \bQ)
\end{align*}
of a class $\gamma \in H^{p} (M, \bQ)$ where $\pi_P$ and $\pi_M$ are projections on $P_n (M, \beta) \times M$.

\begin{definition}
Given $\alpha_j\in \bN$ and $\gamma_j \in H^\ast (M, \bQ)$ for $1 \le j \le r$, the corresponding descendent PT-invariant is
\begin{align*}
    \langle \tp_{\alpha_1 - 1} (\gamma_1) \cdots \tp_{\alpha_r - 1} (\gamma_r) \rangle_{n, \beta} = \int_{[P_n (M, \beta)]^\vir} \prod_{j = 1}^r \tp_{\alpha_j - 1} (\gamma_j).
\end{align*}
We define the following associated partition function
\begin{align}\label{eqn;PT_asso_fun}
    \ZPT \left(M; q \tmid \prod_{j = 1}^r \tp_{\alpha_j - 1} (\gamma_j) \right)_{\beta} = \sum_{n \in \bZ} \left\langle \prod_{j = 1}^r \tp_{\alpha_j-1} (\gamma_j) \right\rangle_{n, \beta} q^n.
\end{align}
\end{definition}

Note that the moduli space $P_n (M, \beta)$ is empty for $n$ sufficiently negative. Therefore \eqref{eqn;PT_asso_fun} is a Laurent series $\bQ\Laurent{q}$ as well.

For the case of primary insertions, Pardon \cite{pardon2023} proved the conjecture of Maulik--Nekrasov--Okounkov--Pandharipande \cite{MNOPI, MNOPII} for all complex threefolds with nef anti-canonical bundle.

\begin{theorem}[{\cite[Theorem 1.6]{pardon2023}}]\label{thm;Pardon}
Let $M$ be a smooth projective threefold with nef anti-canonical bundle. Then, for classes $(\gamma_1, \dots, \gamma_r) \in H^\ast (M, \bQ)^{\oplus r}$, the partition function  $\ZPT \left(M; q \mid \gamma_1  \dots \gamma_{r} \right)_{\beta} \in \bQ (q)$ and 
\[
    (- q)^{- \sC_\beta / 2} \ZPT \left(M  ; q\tmid \gamma_1  \cdots \gamma_{r} \right)_\beta 
     = (- iu)^{\sC_\beta} \ZGWp \left(M; u\tmid \gamma_1 \cdots \gamma_{r} \right)_\beta 
\]
under the variable change $-q = e^{iu}$.
\end{theorem}

To relate descendent GW and PT-invariants, we need the correspondence matrices found by Pandharipande and Pixton \cite{PP14,PP17}. The matrices relating them were predicted in \cite[Conjecture 4]{MNOPII}.

Let $\halpha = (\halpha_1, \dots , \halpha_{\hell})$, with $\halpha_1 \ge \dots \ge \halpha_{\hell} \ge 1$, be a partition of length $\ell (\halpha) \coloneqq \hell$ and size $|\halpha| \coloneqq \sum \halpha_j$. 
Let $\iota_{\Delta} \colon \Delta \to M^{\hell}$ be the inclusion of the small diagonal in the product $M^{\hell}$. For $\gamma \in H^\ast (M, \bQ)$, we write 
\begin{align*}
    \gamma \cdot \Delta \coloneqq \iota_{\Delta \ast} (\gamma) \in H^\ast (M^{\hell}, \bQ).
\end{align*}
Let $\{\theta_j\}$ be a basis of $H^\ast (M, \bQ)$. By K\"unneth formula, we have
\begin{align*}
    \gamma \cdot \Delta = \sum_{j_1, \dots, j_{\hell}} c^{\gamma}_{j_1, \dots, j_{\hell}} \theta_{j_1} \otimes \dots \otimes \theta_{j_{\hell}}.
\end{align*}
The descendent insertion $\tau_{[\halpha]} (\gamma)$ is defined by \cite[(3)]{PP17}
\begin{align*}
    \tau_{[\halpha]} (\gamma) = \sum_{j_1, \dots, j_{\hell}} c^{\gamma}_{j_1, \dots, j_{\hell}} \tau_{\halpha_1 - 1}(\theta_{j_1}) \cdots \tau_{\halpha_{\hell} - 1}(\theta_{j_{\hell}}).
\end{align*}

The key construction in \cite[\S 0.5]{PP14} is a universal correspondence matrix $\tsK$ indexed by partitions $\alpha$ and $\halpha$ of positive size with 
\begin{align*}
    \tsK_{\alpha, \halpha} \in \bQ [\sqrt{- 1}, c_1, c_2, c_3] \Laurent{u}
\end{align*}
and $\tsK_{\alpha, \halpha} = 0$ if $|\alpha| < |\halpha|$. By specializing the formal variables $c_i$ to $c_i (T_M)$, the elements of $\tsK$ act by cup product on $H^\ast (M, \bQ)$ with $\bQ[i]\Laurent{u}$-coefficients.

Let $\alpha = (\alpha_1, \dots, \alpha_{\ell})$ be a partition and $P$ a partition of $\{1, \dots, \ell\}$. For each $S \in P$, a subset of $\{1, \dots, \ell\}$, let $\alpha_S$ be the subpartition consisting of the parts $\alpha_j$ for $j \in S$ and 
\[
    \gamma_S = \prod\nolimits_{j \in S} \gamma_j.
\]

\begin{definition}[\cite{PP14}]\label{inser_def_abs}
For even cohomology classes $\gamma_j \in H^{2 \ast} (M, \bQ)$, let
    \begin{align*}
      \overline{\tau_{\alpha_1 - 1} (\gamma_1) \cdots \tau_{\alpha_{\ell} - 1} (\gamma_{\ell})} = \sum_{\substack{P \text{ set partitions } \\ \text{of } \{1, \dots, \ell \}}} \prod_{S \in P}  \sum_{0 < |\halpha| \leqslant |\alpha_S|} \tau_{[\halpha]} \left(\tsK_{\alpha_S, \halpha} \cdot \gamma_S \right).
    \end{align*}
\end{definition}

\begin{example}
If $\alpha = (1, \ldots,1)$ then $\overline{\tau_{0} (\gamma_1) \cdots \tau_{0} (\gamma_{\ell})} = \tau_{0} (\gamma_1) \cdots \tau_{0} (\gamma_{\ell})$. If $\gamma =\pt$, the class of a point, then 
\[
    \tau_{[\alpha]} (\pt) = \tau_{\alpha_1 - 1} (\pt) \cdots \tau_{\alpha_{\ell} - 1} (\pt).
\]
If $\alpha = (\alpha_1)$, then $\tau_{[\alpha]} (\gamma) = \tau_{\alpha_1 - 1} (\gamma)$. 
\end{example}

\begin{notation}\label{nota;coh}
Given even cohomology classes $\gamma_1, \ldots, \gamma_r \in H^\ast (M, \bQ)$, integers $\alpha_i\in \bN$ for $1\le j\le r$, we set $A = {\prod_{i=1}^r\tau_{\alpha_i-1}(\gamma_i)}$.
\end{notation}

We are now in a position to state the conjectural GW/PT correspondence \cite{PP14}.

\begin{conjecture}\label{GWPTconj}
Let $A$ be as in Notation \ref{nota;coh}. Then:
\begin{enumerate}
    \item\label{conj;rat} The $\ZPT \left(M; q \mid A \right)_{\beta}$ is the Laurent expansion of a rational function in $q$. 

    \item\label{conj;corresp} We have
    \begin{align*}
    (- q)^{- \sC_\beta / 2} \ZPT \left(M  ; q\tmid A \right)_\beta
     = (- iu)^{\sC_\beta} \ZGWp \left(M; u\tmid \oA  \right)_\beta 
    \end{align*}
    under the variable change $- q = e^{iu}$.
\end{enumerate}
\end{conjecture}

Note that the variable change in \eqref{conj;corresp} is well-defined assuming \eqref{conj;rat}, and \eqref{conj;rat} is called the \emph{rationality conjecture}.

%%%%%%%%%%%%%%%%%%%%%%%%%%%%%%%%%%%%%%%%
\subsection{Relative Theories}\label{subsusec;rel_theories}
%%%%%%%%%%%%%%%%%%%%%%%%%%%%%%%%%%%%%%%%

\begin{definition}
Let $D$ be a smooth (connected) divisor on $M$, and let $\cB$ be a basis of $H^\ast (D, \bQ)$. A {\em cohomology weighted partition} $\eta$ with respect to $\cB$ is a set of pairs 
\begin{align*}
    \{(a_1, \delta_1), \dots, (a_{r}, \delta_{r})\}, \quad \mbox{where } \delta_j \in \cB \mbox{ and } a_1 \ge \dots \ge a_{r} \ge 1,
\end{align*}
such that $\vec{\eta} \coloneqq (a_j) \in \bN^{r}$ is a partition of size $|\eta| = \sum a_j$ and length $\ell (\eta) = r$. Its \emph{dual partition} $\eta^\vee$ is the cohomology weighted partition $\{(a_j, \delta_j^\vee)\}_j$ (with respect to the dual basis $\cB^\vee$ of $\cB$). We write $\eta = (1,\delta)^r$ when $a_j = 1$ and $\delta_j = \delta$ for all $j$.
\end{definition}

The automorphism group $\Aut (\eta)$ consists of $\sigma \in \mathfrak{S}_{\ell (\eta)}$ such that $\eta^{\sigma} = \eta$. We set 
\begin{align*}%\label{glufac}
    \fz (\eta) = |{\Aut (\eta)}| \cdot \prod\nolimits_{j = 1}^{\ell (\eta)} a_j.
\end{align*}

\begin{definition}\label{def:relGW}
    Let $A$ be as in Notation \ref{nota;coh} and $\eta$ a cohomology weighted partition. The relative descendent GW-invariant \cite[p.240]{Li02} associated to $A$ and $\eta$ is defined as
    \[
        \langle A \mid  \eta \rangle^\pp_{g, \beta} = \frac{1}{|{\Aut (\eta)}|} \int_{[\oMp_{g, r} (M / D, \beta, \eta)]^{\vir}} \prod_{j = 1}^r \left(\psi_j^{\alpha_j-1} \cup  \ev_j^\ast (\gamma_j) \right) \cup \ev_{D}^\ast (\delta_j).
    \]
    The associated partition function is the Laurent series
    \begin{align}\label{relZGW}
        \ZGWp \left(M / D ; u \mid A \mid \eta\right)_{\beta} 
        = \sum_{g \in \bZ} \left\langle A \mid \eta \right\rangle^\pp_{g, \beta} u^{2 g - 2}
    \end{align}
\end{definition}

In relative PT-theory, we consider the moduli space $P_n (M / D, \beta)$ introduced by Li-Wu \cite{LW15} (cf.\ \cite[\S 3.2]{MNOPII}) which parametrizes stable pairs $(F,s)$ relative to $D$, such that $\chi (F) = n \in \bZ$ and $\ch_2 (F) = \beta$. We have the intersection map
\begin{align}\label{relPTintmap}
    \epsilon \colon P_n (M/D, \beta) \to \Hilb(D, |\eta|)
\end{align}
to the Hilbert scheme of $|\eta| = (D, \beta)$ points of the connected divisor $D$. Fix $d \in \bN$ and let $\eta = \{ (a_j, \delta_j)\}_j$ be a cohomology weighted partition of size $d$ with respect to $\cB$. Let 
\begin{align*}
    C_\eta = \frac{1}{\fz (\eta)} P_{\delta_1} [a_1] \dots P_{\delta_{\ell (\eta)}} [a_{\ell(\eta)}] \cdot \mathbf{1} \in H^{\ast} (\Hilb(D, d ), \bQ),
\end{align*}
see \cite[\S 3.2.2]{MNOPII}. Here $\mathbf{1}$ is the vacuum vector $\vacuum=1\in H^0(\Hilb(D, 0), \bQ)$.
Then $\{C_\eta\}_{|\eta| = d}$ is the Nakajima basis of $H^{\ast} (\Hilb(D,d), \bQ)$ with Poincar\'e pairing 
\begin{align*}
    \int_{\Hilb(D,d)} C_\eta \cup C_\nu = 
    \begin{cases}
        \frac{(- 1)^{d - \ell (\eta)}}{\fz (\eta)} & \text{if }\nu=\eta^\vee, \\
        0 & \text{otherwise.}
    \end{cases} 
\end{align*}

\begin{definition}
    With notation as in Definition~\ref{def:relGW}, the relative descendent PT-invariant associated to $A$ and $\eta$ is
    \begin{align*}
        \langle A \mid \eta \rangle_{n, \beta} = \int_{[P_n (M / D, \beta)]^{\vir}} \left(\prod_{i = 1}^r \tp_{\alpha_i - 1} (\gamma_i) \right) \cup \epsilon^\ast (C_{\eta}).
    \end{align*}
    The associated partition function is the Laurent series
    \begin{align}\label{relZPT}
         \ZPT \left(M / D ; q\mid A \mid \eta \right)_{\beta} 
         = \sum_{n \in \bZ} \left\langle A \mid \eta \right\rangle_{n, \beta} q^n.
    \end{align}
\end{definition}

Now, we can state the conjectural relative descendent GW/PT correspondence \cite{MNOPII, PP14, PP17}.

\begin{conjecture}\label{relGWPTconj}
With notation as in Definition~\ref{def:relGW}, we have:
\begin{enumerate}
    \item The $\ZPT \left(M / D; q\tmid A \tmid \eta  \right)_{\beta}$ is the Laurent expansion of a rational function in $q$.

    \item Under the variable change $e^{i u} = -q$,
    \begin{align*}
    (-q)^{- \sC_\beta^M / 2} \ZPT \left(M / D; q \tmid A \tmid \eta\right)_\beta
    = (- i u)^{\sC_\beta^M + \ell (\eta) - |\eta|} \ZGWp \left(M / D; u \tmid \oA \tmid \eta \right)_\beta.
    \end{align*}
\end{enumerate}
\end{conjecture}

To study Conjectures \ref{GWPTconj} and \ref{relGWPTconj} for $Y_d  = \bP_{S_d} (K_{S_d} \oplus \cO)$ over a nontoric surface $S_d$, we need the following proposition.

\begin{proposition}\label{prop;def_to_toric}
For each $1 \leq d \leq 8$, there exists a toric surface $S^0_d$ such that it is deformation equivalent to $S_d$. 
\end{proposition}

\begin{proof}
For $d \in \{6, 7, 8\}$, $S_d$ is already a toric surface. Hence, we assume that $d \le 5$. Consider the points 
\begin{align*}
    &\sigma_1(t) = [1: t^5: 0], \qquad \sigma_2(t) = [t^5: 1: 0], \qquad \sigma_3(t) = [0: 0: 1]\\
    &\sigma_4(t) = [1: t^2: t], \qquad \sigma_5(t) = [t^2: 1: t], \qquad \sigma_6(t) = [t: t^2: 1] \\
    &\sigma_7(t) = [1: t^3: 2t], \qquad \sigma_8(t) = [t^4: 1: 2t]. 
\end{align*}
It is checked that for generic $t \neq 0$, the points $\{\sigma_i(t)\}_{1\le i \le 8} \subseteq \bP^2$ lie in a general position, i.e., 
\begin{inparaenum}[(i)]
    \item no three are collinear; 
    \item no six on the same conic; 
    \item no eight on a cubic with a double point at one of them. 
\end{inparaenum} 
Consider the deformation family defined by a sequence of blow-ups
\[\cS = \cS''' \xrightarrow{\ \pi''\ } \cS'' \xrightarrow{\ \pi'\ } \cS' \xrightarrow{\ \pi\ } \bP^2 \times \bC \mathop{\longrightarrow}^t \bC, \]
where 
\begin{itemize}
    \item $\pi$ is the blow-up along the graphs of $\sigma_1$, $\sigma_2$, $\sigma_3$; 
    \item $\pi'$ is the blow-up along the strict transform of the graphs of $\sigma_4$, $\sigma_5$, $\sigma_6$; 
    \item $\pi''$ is the blow-up along the strict transform of the graphs of $\sigma_7$, $\sigma_8$. 
\end{itemize}
It follows that the central fiber of $\cS$ is the toric surface defined by the rays
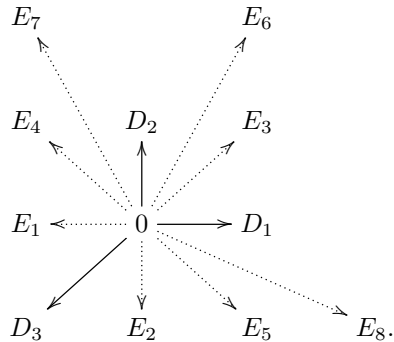
\[\mbox{\footnotesize
    \xymatrix{
    E_7 & & E_6 \\ 
    E_4 & D_2 & E_3 \\ 
    E_1 & 0 \ar[r] \ar[u] \ar[dl] \ar@{.>}[l] \ar@{.>}[lu] \ar@{.>}[luu] \ar@{.>}[d] \ar@{.>}[rd] \ar@{.>}[rrd] \ar@{.>}[ru] \ar@{.>}[ruu] & D_1 \\ 
    D_3 & E_2 & E_5 & E_8.
    }}
\]
Also, if we only blow-up the strict transforms of the graphs $\sigma_1$, $\dots$, $\sigma_k$, then the general fiber is a del Pezzo surface of degree $d = 9 - k$ and the central fiber is the toric surface $S_{9-k}^0$ defined by the rays $D_1$, $D_2$, $D_3$, $E_1$, $\dots$, $E_k$. 
\end{proof}

\begin{remark}
For $d \geq 3$, one can compare Proposition \ref{prop;def_to_toric} with \cite[Proposition 4.1]{KM09}.
\end{remark}

%%%%%%%%%%%%%%%%%%%%%%%%%%%%%%%%%%%%%%%%
\subsection{Correspondence}
%%%%%%%%%%%%%%%%%%%%%%%%%%%%%%%%%%%%%%%%

First, we apply Proposition \ref{prop;def_to_toric} to prove the descendent GW/PT correspondence for certain threefolds.

\begin{theorem}\label{thm;GWPTcorrp_def_toric}
Conjecture \ref{GWPTconj} holds for $M$ being one the following threefolds:
\begin{enumerate}
    \item\label{thm;GWPTcorrp_def_toric_1} the projective bundle $Y_d = \bP_{S_d} (K_{S_d} \oplus \cO)$; 
    \item\label{thm;GWPTcorrp_def_toric_2} a smooth Fano threefolds with Picard number $\rho \geq 6$.
\end{enumerate}
\end{theorem}

\begin{proof}
By Proposition~\ref{prop;def_to_toric}, we can deform $S_d$ to a toric surface $S_d^0$. Thus $Y_d = \bP_{S_d} (K_{S_d} \oplus \cO)$ and $Y_d^0 = \bP_{S_d^0} (K_{S_d^0} \oplus \cO)$ are also deformation equivalent. On the other hand, if $M$ satisfies \eqref{thm;GWPTcorrp_def_toric_2}, then $M = \bP^1 \times S_{\rho - 2}$ by the Mori--Mukai classification \cite{MM82}. Then the theorem follows from \cite[Theorem 7]{PP14} and the deformation invariance of GW and PT-invariants.
\end{proof}

\begin{remark}\label{rmk;exa}
The Fano threefolds $\bP^1 \times S_{1}$ and $\bP^1 \times S_{5}$ are not complete intersections in product of projective spaces, and thus they are not treated in \cite{PP17}.
\end{remark}

We also establish the relative descendent GW/PT correspondence for $Y_d$, which will be used in Theorem \ref{thm;GWPT}.

\begin{theorem}\label{thm;rel_loc_GWPTdPs}
For $1 \leq d \leq 8$, let $E$ and $H$ be the zero section $\bP_{S_d}(\cO)$ and the infinity section $\bP_{S_d}(K_{S_d})$ of $Y_d$ respectively. Then Conjecture~\ref{relGWPTconj} holds for $Y_d/E$ and $Y_d/H$.
\end{theorem}

\begin{proof}
As before, $S_d$ and $S_d^0$ are deformation equivalent by Proposition~\ref{prop;def_to_toric}. Let $E^0$ (resp.~$H^0$) be the zero section $\bP_{S_d^0}(\cO)$ (resp.~infinity section $\bP_{S_d^0}(K_{S_d^0})$) of $Y_d^0 = \bP_{S_d^0} (K_{S_d^0} \oplus \cO)$. It follows from \cite[Theorem~2]{PP17} that Conjecture~\ref{relGWPTconj} holds for the pairs $Y_d^0/E^0$ and $Y_d^0/H^0$. Since both relative GW-invariants and relative PT-invariants are invariant under deformations of pairs, Conjecture~\ref{relGWPTconj} also holds for $Y_d/E$ and $Y_d/H$. 
\end{proof}

To use the double point degenerations \eqref{eqn;ss_Y} and \eqref{eqn;ss_X}, we need the following degeneration formulas for GW and PT-invariants. To save notations, we denote these two degenerations by $\cW \to \Delta$. It has a smooth fiber $M \coloneqq \cW_t$ ($t \neq 0$), a special fiber $\cW_0 = M_0 \cup_D M_\infty$ and $D = M_0 \cap M_\infty$ a smooth divisor. Let $\iota \colon M \hookrightarrow \cW$, $\iota_0 \colon M_0 \hookrightarrow \cW$ and $\iota_\infty \colon M_\infty \hookrightarrow \cW$ be the inclusion maps.

\begin{theorem}\label{thm:deg-gw}
    Suppose that $\gamma_1$, $\dots$, $\gamma_r$ are even cohomology classes on the total space $\cW$, and let $A = \tau_{\alpha_1-1}(\gamma_1)\dots\tau_{\alpha_r-1}(\gamma_r)$. 
    For a nonzero class $\beta^\prime \in \NE (\cW)$, we have
    \begin{align*}
        &\sum_{\substack{
        \beta \in \NE (M) \\ 
        \iota_\ast \beta = \beta^\prime}} \ZGWp \left(M ;u\mid \overline{A} \right)_{\beta} \\
        &\quad =
        \sum \fz(\eta)  u^{2 \ell(\eta)} \ZGWp \left(M_0 / D; u \mid \overline{A_{I_0}} \mid \eta \right)_{\beta_0} \ZGWp \left(M_\infty / D; u \mid \overline{A_{I_\infty}} \mid \eta^{\vee} \right)_{\beta_\infty}, 
    \end{align*}
    where the summation on the second line runs over
    \begin{enumerate}[(a)]
        \item splittings 
        \begin{align}\label{eq:curve-splitting-general}
            \iota_{0 \ast} \beta_0 + \iota_{\infty \ast} \beta_\infty = \beta^\prime = \iota_*\beta
        \end{align} 
        such that $(D, \beta_0) = (D, \beta_\infty)$,
        \item partitions $I_0 \sqcup I_\infty =\{1,2, \dots, r\}$, and 
        \item cohomology weighted partition $\eta$ such that $|\eta|= (D, \beta_0)$ with respect to a fixed basis of $H(D, \bQ)$. 
    \end{enumerate} 
\end{theorem}
See for example \cite{Li02} and \cite[p.403]{PP17}. Similarly, the formula without bars, namely without applying the universal transformation to descendent insertions, also holds.

\begin{theorem}\label{thm:deg-pt}
    With notation as in Theorem~\ref{thm:deg-gw}, we have
    \begin{align*}
        &\sum_{\substack{\beta \in \NE (M) \\ \iota_\ast \beta = \beta^\prime}} \ZPT \left(M; q \mid A \right)_{\beta} \\
        &\quad =\sum (-1)^{\ell (\eta)} \fz(\eta) (-q)^{- |\eta|}\ZPT \left(M_0 / D; q \mid 
        A_{I_0} \mid \eta\right)_{\beta_0} \ZPT \left(M_\infty/ D; q \mid A_{I_\infty} \mid \eta^{\vee} \right)_{\beta_\infty}, 
    \end{align*}
    where the summation on the second line runs over the same index set in Theorem \ref{thm:deg-gw}. 
\end{theorem}

See for example \cite{LW15, MNOPII}, \cite[p.2761]{PP14} and \cite[Theorem 6.12]{Lin23}.

Now, we are ready to prove the descendent GW/PT correspondence for del Pezzo transitions.

\begin{theorem}\label{thm;GWPT}
Let $\typeII$ be a del Pezzo transition with the smoothing $\fX \to \Delta$. Suppose that $\beta\in \NE(X)$ is a nonzero class and $\alpha=(\alpha_1,\dots, \alpha_r)$ a fixed partition. Assume $\gamma_i\in H^{\ast}(\fX, \bQ)$, $i=1$, $\dots$, $r$, are fixed even cohomology classes and if $\gamma_i\in H^0(\mathfrak{X}, \bQ)$, then $\alpha_i=1$. 
\begin{enumerate}[(a)]
    \item\label{thm:GWPT_1} If Conjecture~\ref{GWPTconj} \eqref{conj;rat} holds for $\Xres$, then it holds for $X$ and descendent insertions 
        \begin{align}\label{eqn;ins-total}
            \gamma_{i|X},\quad i=1,\dots ,r.
        \end{align}
        \item\label{thm:GWPT_2} Set $A = {\prod_{i=1}^r\tau_{\alpha_i-1}(\gamma_{i}|_X)}$. If furthermore Conjecture~\ref{GWPTconj} \eqref{conj;corresp} holds for $\Xres$, then it holds for $\Xsm$ with descendent insertions \eqref{eqn;ins-total}, i.e.,
        \begin{align*}
        (- q)^{- \sC_\beta / 2} \ZPT \left(X; q \tmid A \right)_\beta 
        = (- iu)^{\sC_\beta} \ZGWp \left(X; u\tmid \overline{A}  \right)_\beta. 
        \end{align*}
    \end{enumerate}
\end{theorem}

\begin{proof}
By the string equation, we may assume that $\gamma_i \in H^{>0}(\fX, \bQ)$ for all $i$. Let $d = E^3$ be the degree of the transition $\Xres \searrow \Xsm$. Set $S_\loc = S_d$, $X_\loc = X_d$, and $Y_\loc = Y_d$. We denote by $h = c_1(\cO_{X_\loc}(1))$, and by $E$ (resp.~$H$) the zero section $\bP_{S_\loc}(\cO)$ (resp.~the infinity section $\bP_{S_\loc}(K_{S_\loc})$) of 
$\pi \colon Y_\loc \to S_\loc$.

Applying the degeneration formulas, Theorems \ref{thm:deg-gw} and \ref{thm:deg-pt}, to the double point degeneration \eqref{eqn;ss_X}, and a virtual dimension counting, we get  
\begin{align*}
    \ZGWp\left(X;u \mid \overline{A}\right)_{\beta} &=
    \sum_\rho
    \rho! u^{2\rho} \ZGWp\left(Y / E; u \mid \overline{\phi^*A} \mid (1, 1)^\rho\right)_{\beta_0}
    \ZGWp\left(\Xloc / h ; u  \mid \varnothing \mid (1, \pt)^\rho \right)_{\beta_\infty}, \\
    \ZPT\left(X;q \mid A\right)_{\beta} &=
    \sum_\rho
    \rho! q^{-\rho} \ZPT\left(Y / E; u \mid \phi^*A \mid (1, 1)^\rho\right)_{\beta_0}
    \ZPT\left(\Xloc / h ; q\mid \varnothing \mid (1, \pt)^\rho \right)_{\beta_\infty} 
\end{align*}
where $(E, \beta_0) = (h, \beta_\infty) = \rho$ and $\phi_*\beta_0 = \beta$.

For $r \ge 0$ and a cohomology weighted partition $\eta$, we define
\begin{align*}
    \ZGWp\left(\Xloc/ h; u \mid A \mid \eta\right)_{r} &\coloneqq \sum_{(h, \beta_\loc) = r}\ZGWp\left(\Xloc/ h; u \mid A \mid \eta\right)_{\beta_\loc}, \\
    \ZPT\left(\Xloc / h; q \mid A \mid \eta\right)_{r} &\coloneqq \sum_{(h, \beta_\loc) = r}\ZPT\left(\Xloc / h; q \mid A \mid \eta\right)_{\beta_\loc}. 
\end{align*}
Then it suffices to show that the conjecture holds for $Y/E$ with any insertion and for $\Xloc/h$ with trivial insertion, i.e., for $(E, \tilde{\beta}) = r$, 
\begin{align}
    (-iu)^{\sC_\beta} \ZGWp\left(Y / E; u \mid \overline{\phi^*A} \mid (1, 1)^r\right)_{\tilde{\beta}} &=  (-q)^{\sC_\beta/2}\ZPT\left(Y / E; q \mid \phi^*A \mid (1, 1)^r\right)_{\tilde{\beta}}, \label{eq:P(Y/E)}\\
    (-iu)^{2r} \ZGWp\left(\Xloc/ h; u \mid \varnothing \mid (1, \pt)^r\right)_{r} &=  (-q)^r \ZPT\left(\Xloc / h; q \mid \varnothing \mid (1, \pt)^r\right)_{r}, \label{eq:P(Xloc/h)}
\end{align}
and are rational in $q$.

For \eqref{eq:P(Y/E)}, we apply the degeneration formulas to the double point degeneration \eqref{eqn;ss_Y}:
\begin{align*}
    \ZGWp\left(Y;u \mid \overline{\phi^*A}\right)_{\tilde{\beta}} &=
    \sum_\rho
    \rho! u^{2\rho} \ZGWp\left(Y / E; u \mid \overline{\phi^*A} \mid (1, 1)^\rho\right)_{\tilde{\beta}_0}
    \ZGWp\left(\Yloc / H ; u  \mid \varnothing \mid (1, \pt)^\rho \right)_{\tilde{\beta}_\infty}, \\
    \ZPT\left(Y;q \mid \phi^*A\right)_{\tilde{\beta}} &=
    \sum_\rho
    \rho! q^{-\rho} \ZPT\left(Y / E; q \mid \phi^*A \mid (1, 1)^\rho\right)_{\tilde{\beta}_0}
    \ZPT\left(\Yloc / H ; q\mid \varnothing \mid (1, \pt)^\rho \right)_{\tilde{\beta}_\infty}, 
\end{align*}
where $(E, \tilde{\beta}_0) = (H, \tilde{\beta}_\infty) = \rho$ and $\tilde{\beta}_0 + \pi_*\tilde{\beta}_\infty = \tilde{\beta}$. Let $\ff$ be the fiber class of $\pi \colon Y_\loc \to S_\loc$. Note that when $\tilde{\beta}_0 = \tilde{\beta}$, we must have $\rho = r$ and $\beta_\infty = r \ff$. 

By induction on the order of $\tilde{\beta}_0$, we only need to show that
\begin{enumerate}[(i)]
    \item\label{eqn;induc_1} the conjecture holds for $\Yloc/H$; 
    \item\label{eqn;induc_2} the terms $\ZGWp\left(\Yloc/H; u \mid \varnothing \mid (1, \pt)^r\right)_{r \ff}$, $\ZPT\left(\Yloc/H; q \mid \varnothing \mid (1, \pt)^r\right)_{r \ff}$ are not zero. 
\end{enumerate}
The \eqref{eqn;induc_1} follows from Theorem~\ref{thm;rel_loc_GWPTdPs}. For \eqref{eqn;induc_2}, we apply \cite[Proposition~6]{PP13}, which gives 
\[\ZPT\left(\Yloc/H; q \mid \varnothing \mid (1, \pt)^r\right)_{r \ff} = \int_{\Hilb(S, r)} C_{(1, \pt)^r} \neq 0,\]
and, by \eqref{eqn;induc_1}, 
\[(-iu)^{2r}\ZGWp\left(\Yloc/H; u \mid \varnothing \mid (1, \pt)^r\right)_{r \ff} = (-q)^r\ZPT\left(\Yloc/H; q \mid \varnothing \mid (1, \pt)^r\right)_{r \ff} \neq 0. \]

For \eqref{eq:P(Xloc/h)}, we apply the degeneration formula to $\Xloc \rightsquigarrow \Yloc \cup_{E} \Xloc$: if we let 
\begin{align*}
    \ZGWp\left(\Xloc; u \mid A\right)_{r} &\coloneqq \sum_{(h, \beta_\loc) = r}\ZGWp\left(\Xloc; u \mid A\right)_{\beta_\loc}, \\
    \ZPT\left(\Xloc; q \mid A\right)_{r} &\coloneqq \sum_{(h, \beta_\loc) = r}\ZPT\left(\Xloc; q \mid A\right)_{\beta_\loc},
\end{align*}
then 
\begin{align*}
    \ZGWp\left(\Xloc ; u \mid \pt^r\right)_{r} &=
    \sum_{\rho}
    \rho! u^{2\rho} \ZGWp\left(\Yloc / E; u \mid \pt^r \mid (1, 1)^\rho\right)_{\beta_0}
    \ZGWp\left(\Xloc / h ; u  \mid \varnothing \mid (1, \pt)^\rho \right)_{\rho}, \\
    \ZPT\left(\Xloc ; q \mid \pt^r\right)_{r} &=
    \sum_{\rho}
    \rho! q^{-\rho} \ZPT\left(\Yloc / E; q \mid \pt^r \mid (1, 1)^\rho\right)_{\beta_0}
    \ZPT\left(\Xloc / h ; q\mid \varnothing \mid (1, \pt)^\rho \right)_{\rho}, 
\end{align*}
where $(E, \beta_0) = \rho$ and $(h, \phi_*\beta_0) = r$. Since $H - E = \pi^*c_1(\cO_{S_d}(1))$ is nef on $Y_\loc$, the condition implies that 
\[0 \le (H - E, \beta_0) = (h, \phi_*\beta_0) - (E, \beta_0) = r - \rho. \]
By induction on $r$, it suffices to prove
\begin{enumerate}[\ (1)]
    \item\label{eqn;induc_loc_1} the conjecture holds for $\Xloc$ with no $\psi$-classes; 
    \item\label{eqn;induc_loc_2} the conjecture holds for $\Yloc/E$; 
    \item\label{eqn;induc_loc_3} the terms $\ZGWp\left(\Yloc/E; u \mid \pt^r \mid (1, 1)^r\right)_{r\ff}$, $\ZPT\left(\Yloc/E; q \mid \pt^r \mid (1, 1)^r\right)_{r \ff}$ are not zero. 
\end{enumerate}
The \eqref{eqn;induc_loc_1} and \eqref{eqn;induc_loc_2} follow from Theorems \ref{thm;Pardon} and \ref{thm;rel_loc_GWPTdPs} respectively. For \eqref{eqn;induc_loc_3}, we apply \cite[Proposition~6]{PP13} again, which gives
\[\ZPT\left(\Yloc/E; q \mid \pt^r \mid (1, 1)^r\right)_{r \ff} = \int_{\Hilb(S, r)} (\tau_0(\pt))^r \cdot C_{(1,1)^r} \neq 0, \]
and, by \eqref{eqn;induc_loc_2}, 
\[(-iu)^{2r}\ZGWp\left(\Yloc/E; u \mid \pt^r \mid (1, 1)^r\right)_{r\ff} = (-q)^r\ZPT\left(\Yloc/E; q \mid \pt^r \mid (1, 1)^r\right)_{r \ff} \neq 0. \]
This completes the proof.
\end{proof}

\begin{remark}\label{rmk;ext_classes}
Given a del Pezzo transition $\typeIIno$, we assume that $h^i (\cO_{\Xsm}) = h^i(\cO_{\Xres}) = 0$ for $i = 1, 2$. Then $\rho(\Xsing) = \rho(\Xsm)$ (cf.\ \cite[Proposition 3.1]{Kapustka209I}) and thus one may infer that the restriction map $H^{\mathrm{ev}}(\fX, \bQ) \cong H^{\mathrm{ev}}(\Xsing, \bQ) \to H^{\mathrm{ev}}(\Xsm, \bQ)$ is surjective (cf.\ \cite[Proposition 1.13]{LWW25}).
\end{remark}

In the proof of Theorem \ref{thm;GWPT}, we have used Pardon's result (Theorem \ref{thm;Pardon}) on $X_d$ for primary insertions. Applying Theorem \ref{thm;GWPT}, we can further prove that Conjecture \ref{GWPTconj} holds for $X_d$ with stationary descendent insertions.

\begin{corollary}\label{cor;GWPT_dP3}
Every smooth del Pezzo threefold $X_d$ satisfies the descendent $\GW/\PT$ correspondence with stationary descendent insertions, i.e., all descendent insertions are even classes of positive degree. 
\end{corollary}

\begin{proof}
For each $1 \leq d \leq 8$, consider the del Pezzo trasnition $\Xres_{d} \searrow \Xsm_{d}$ construced in Example~\ref{ex:dPtrans}. By Theorem~\ref{thm;GWPTcorrp_def_toric}, $\Xres_{d}$ satisfies satisfies Conjecture \ref{GWPTconj}. Since $h^i (\cO_{\Xsm_d}) = h^i(\cO_{\Xres_d}) = 0$ for $i = 1, 2$, every even cohomology class of $\Xsm$ has a lifting in $H^{\mathrm{ev}}(\Xsing, \bQ) \cong H^{\mathrm{ev}}(\fX, \bQ)$, see Remark \ref{rmk;ext_classes}. Therefore $X_d$ satisfies the descendent $\GW/\PT$ correspondence with stationary descendent insertions by Theorem~\ref{thm;GWPT}.
\end{proof}

\begin{remark}\label{rmk;exaX1}
The $X_1$ is a new case which is not contained in \cite{LW25, PP17}. For $d \in \{2,5\}$, $X_d$ is treated in \cite[Corollaries 4.2 \& 4.5]{LW25} and other cases are contained in \cite{PP14, PP17}.
\end{remark}

\end{document}